\documentclass[11pt]{amsart}
\usepackage{geometry}                
\usepackage{graphicx}
\usepackage{amssymb}
\usepackage{epstopdf}
\newcommand{\R}{\mathbb{R}}
\newcommand{\N}{\mathbb{N}}

\title[A family of Gaussian processes]{A natural 4-parameter family of covariance functions for stationary Gaussian processes}
\author{R.S.MacKay \and N.E.Phillips}
\email{R.S.MacKay@warwick.ac.uk \and nicholas.phillips@epfl.ch}
\address{Mathematics Institute, University of Warwick, Coventry CV4 7AL, U.K. \and  
Institute of Bioengineering, School of Life Sciences, Ecole Polytechnique FŽdŽrale de Lausanne, Lausanne CH-1015, Switzerland}

\date{\today}                                           


\begin{document}

\begin{abstract}
A four-parameter family of covariance functions for stationary Gaussian processes is presented.  We call it 2Dsys.  It corresponds to the general solution of an autonomous second-order linear stochastic differential equation, thus arises naturally from modelling.  It covers underdamped and overdamped systems, so it is proposed to use this family when one wishes to decide if a time-series corresponds to stochastically forced damped oscillations or a stochastically forced overdamped system.
\end{abstract}

\maketitle

\section{Introduction}
Gaussian processes (GP) are a flexible class of probability distributions over functions, useful for data analysis and inference.  
A GP on a space $T$ (e.g.~$T=\R$, representing time) is specified by a mean function $M: T \to \R$ and a covariance function $C: T\times T \to \R$.  
It gives a probability distribution on functions $f: T \to \R$, such that for any finite set of points $t_i \in T$ with $i = 1,\ldots,n$,
the values $f_i = f(t_i) \in \R$ are jointly Gaussian distributed with mean vector $m_i = M(t_i)$ and covariance matrix $c_{ij} = C(t_i,t_j)$.  

The covariance function is required to be symmetric and positive-definite.  The latter is defined most simply as the matrix $c$ being positive-definite for all finite sets of observation points $t_i$.

GPs often come in families, labelled by parameters (in much literature known as hyperparameters).  Lists of commonly used families are given in \cite{L+,RW}, where it is also explained how they may be combined. 

One example is the Ornstein-Uhlenbeck (OU) process on $\R$.  It can be defined by $M(t)=0$ and 
\begin{equation}
C(s,t) = a^2 \exp{(-\mu |t-s|)},
\end{equation}
with parameters $a,\mu > 0$ (its covariance function is also known as Mat\'ern-$\frac12$).  It is an example of a stationary process, meaning that the probability distribution is invariant under translation in $\R$; equivalently for a GP, $M$ is constant and $C(s,t)$ is a function of $t-s$ (and so of $|t-s|$ by symmetry of covariance functions).

The original description of the OU process is as the solution of the 1D linear stochastic system
\begin{equation}
\dot{x} = -\mu x + \xi
\end{equation}
started at time $t=-\infty$, with over-dot denoting $d/dt$ and $\xi$ Gaussian white noise of zero mean and covariance 
\begin{equation}
\langle \xi(s) \xi(t) \rangle = K \delta(t-s), K>0.
\end{equation}  
One obtains the relation $a^2 = K/2\mu$.
Indeed, linear time-invariant stochastic differential equations generate a large range of stationary GPs.  

Our aim here is to introduce a 4-parameter family of stationary covariance functions for GPs on $\R$.  It arises naturally from second-order stable autonomous stochastic linear systems.  We believe it will have many uses.  In particular, we propose its use for deciding if a noisy signal corresponds to an underdamped or overdamped system.

\section{2D linear stochastic systems}

The general 2D linear autonomous stochastic system can be written as
\begin{equation}
\dot{x} = \left[ \begin{array}{cc} -A & B \\ C & -D \end{array} \right] x + \xi,
\label{eq:system}
\end{equation}
with $\xi$ being a 2D Gaussian white noise of zero-mean and covariance 
$$\langle \xi_i(s) \xi_j(t) \rangle = K_{ij} \delta(t-s)$$ 
with $K$ positive semi-definite (psd), meaning that 
\begin{equation}
u^T K u \ge 0\ \forall u \in \R^2.
\label{eq:psd}
\end{equation}

Assume stability: 
\begin{eqnarray}
A+D &>& 0 \label{eq:tr}\\
AD &>& BC. \label{eq:det}
\end{eqnarray}
The latter can be written as $\det > 0$ with 
\begin{equation}
\det = AD-BC.
\end{equation}

Then the vector response $x$ can be written as the convolution of the matrix impulse response $h$ with the noise vector $\xi$:
\begin{equation}
x(t) = \int_{-\infty}^t h(t-\tau) \xi(\tau)\ d\tau,
\end{equation}
where $h$ is the matrix solution of 
\begin{equation}
\dot{h} = \left[\begin{array}{cc} -A & B \\ C & -D \end{array} \right] h
\end{equation}
for $t\ge 0$ with $h(0)=I$.

Since $\xi$ is assumed Gaussian and $x$ is linear in $\xi$, then $x$ is a Gaussian process.
The mean of $x(t)$ is zero and the covariance matrix $C(t)$ of $x$ can be computed for $t>0$ by
$$C_{ij}(t) = \langle x_i(0)x_j(t) \rangle = \int_0^\infty h_{ik}(s) K_{kl} h_{jl}(s+t)\ ds,$$
so
\begin{equation}
C(t) = S h^T (t),
\label{eq:C}
\end{equation}
with
\begin{equation}
S = \int_0^\infty h(s) K h^T(s)\ ds.
\end{equation}
Similarly, for $t<0$, $$C(t)=h(-t)S.$$

Calculation of the 11-component of the covariance function for $x$ gives
\begin{equation}
C_{11}(t) = e^{-\sigma|t|} \left( S_{11} \cosh\sqrt{\Delta}t + \sigma J \frac{\sinh\sqrt{\Delta}|t|}{\sqrt{\Delta}} \right),
\label{eq:C11}
\end{equation}
where 
\begin{eqnarray}
\sigma &=& (A+D)/2 \\
\Delta &=& (A-D)^2/4+BC = \sigma^2-\det \label{eq:Delta} \\
S_{11} &=& (D^2 K_{11}+2BDK_{12}+B^2K_{22}) + K_{11} \det  \label{eq:S11} \\
J &=&(D^2 K_{11}+2BDK_{12}+B^2K_{22}) - K_{11} \det . \label{eq:J}
\end{eqnarray}

If $\Delta<0$ then 
\begin{eqnarray}
\cosh\sqrt{\Delta} t &=& \cos \sqrt{-\Delta} t \\
\frac{\sinh\sqrt{\Delta}|t|}{\sqrt{\Delta}} &=& \frac{\sin\sqrt{-\Delta}|t|}{\sqrt{-\Delta}}, \nonumber 
\end{eqnarray}
which some computer programs will do automatically and others might need some help to achieve.  Both functions are even in $\sqrt{\Delta}$, so the square root causes no singularity nor imaginary values.
For $\Delta=0$, we obtain the limiting case
\begin{equation}
C_{11}(t) = e^{-\sigma |t|} \left(S_{11} + \sigma J |t| \right).
\end{equation}
which some computer programs do not manage to obtain automatically, so one has to assist them explicitly.

The only constraints on the parameters of $C_{11}$ are $\sigma > 0$, $\Delta < \sigma^2$, $S_{11} > 0$ and $|J| \le S_{11}$.
These constraints follow in turn from the first stability condition (\ref{eq:tr}), (\ref{eq:Delta}) and the second stability condition (\ref{eq:det}), (\ref{eq:S11}) via (\ref{eq:det}) and $K$ psd (take $u = (D,B)$ and $(1,0)$ in (\ref{eq:psd}), and exclude the case $S_{11}=0$ because then $C_{11}(t)=0$ is trivial), and the same two conditions again applied to (\ref{eq:S11}) and (\ref{eq:J}).

Thus (\ref{eq:C11}) gives a 4-parameter family of covariance functions, parametrised by $\sigma>0$, $\Delta < \sigma^2$, $S_{11}>0$ and $J$ in $|J| \le S_{11}$.

It is convenient to reparametrise the family by writing 
\begin{eqnarray}
\sigma &=& \exp(h), \\
\Delta &= &\sigma^2(1-\exp(s)), \\ 
S_{11} &=& \exp(2k), \\
J&=&S_{11} j \mbox{ with } j \in [-1,+1],
\end{eqnarray}
taking new parameters $h, s, k$ and $j$.  
Thus (dropping the subscript $11$ now) the family can be written as
\begin{equation}
C(t) = e^{2k} \exp(-e^h |t|) \left( \cosh(\sqrt{1-e^s}\ e^h t) + j \ \frac{\sinh(\sqrt{1-e^s}\ e^h |t|)}{\sqrt{1-e^s}} \right).
\end{equation}

Uniform priors on the new parameters are a plausible starting point, though knowledge about the system under study should be incorporated.  
For example, if the noise covariance matrix $K$ and the matrix in (\ref{eq:system}) have no special structure then from (\ref{eq:S11},\ref{eq:J}), $j=+1$ is a codimension-1 event ($K_{11} \det = 0$) and $j=-1$ is a codimension-2 event ($K_{12}^2=K_{11}K_{22}$ and $B K_{11}+D K_{12}=0$), so the prior density should have a positive limit as $j \to +1$ and go proportional to $1+j$ as $j \to -1$.
On the other hand, if one knows that $K_{11}=0$ then $j$ is forced to be $+1$.

It may be convenient for optimisation of parameter fits to replace $j$ by an unconstrained parameter.  For example, one could write $j = \sin \pi p/2$ with $p$ unconstrained, and regard $p = n+\frac12 \pm \alpha$ as equivalent for any $n \in \N$ and choice of $\pm$. 
An alternative is to write $j = \tanh r$ with $r$ unconstrained, but this pushes the two permitted cases $j = \pm 1$ off to $r = \pm \infty$, which may or may not be desirable, depending on the context; this can also be written as $j = \frac{1-z}{1+z}$ with $z = e^{-2r}>0$.  Other alternatives of a similar nature are $j = \frac{2}{\pi} \arctan z$ or $j = \frac{2z}{1+\sqrt{1+4z^2}}$.

The parameters $k$ and $h$ just determine scales in $x_1$ and $t$ respectively.  The parameter $s$ determines the separation of the eigenvalues: $s=0$ corresponds to a double eigenvalue (critical damping), $s>0$ to underdamped motion, and $s<0$ to overdamped.  The parameter $j$ determines to what extent $x_1$ is noise-driven directly versus indirectly via $x_2$: $j=+1$ corresponds to no direct noise-driving, $j=-1$ to no indirect noise-driving.  

Related to the parameter $j$, note that the derivative $C'(0+)$ is zero if and only if $j=+1$.  Using \cite{A}, this implies that for $j=+1$ almost every sample is differentiable, whereas for $j\ne +1$ almost every sample is not differentiable.  So $1-j$ measures the extent of non-differentiability in the process.  This is a different type of parameter from $\nu$ in the Mat\'ern covariance functions, which regulates the degree of smoothness, generating samples with $\bar{\nu}-1$ derivatives, where $\bar{\nu}$ denotes the least integer not less than $\nu$.

\begin{figure}[htbp] 
   \centering
   \includegraphics[width=1.12in]{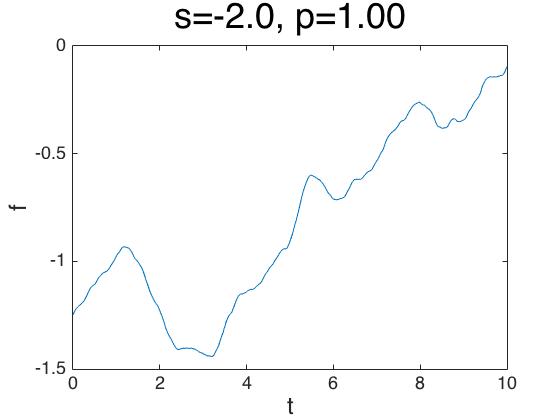} \includegraphics[width=1.12in]{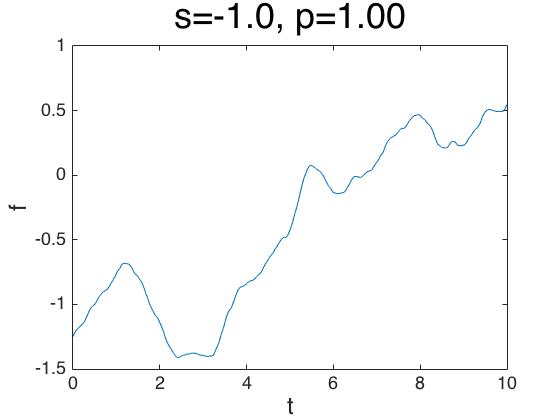} \includegraphics[width=1.12in]{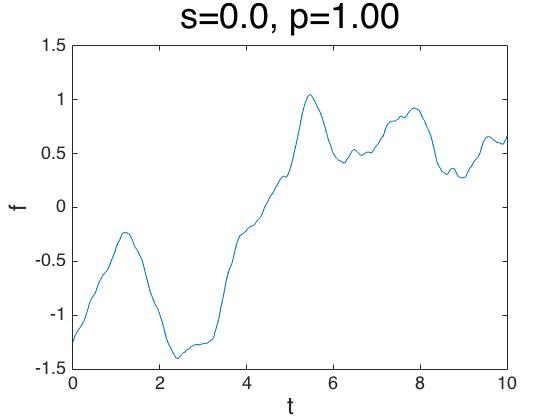} \includegraphics[width=1.12in]{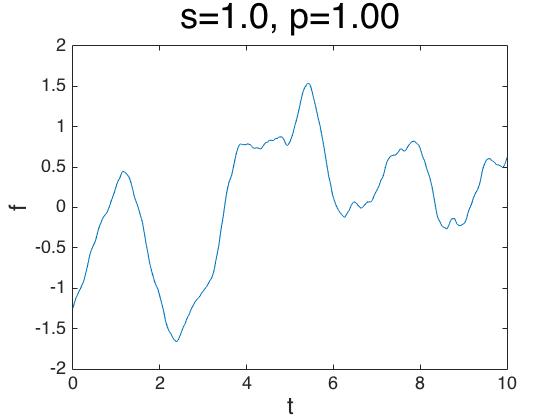} \includegraphics[width=1.12in]{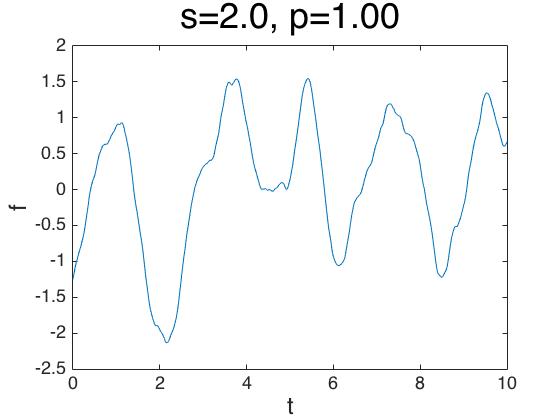} \\
   \includegraphics[width=1.12in]{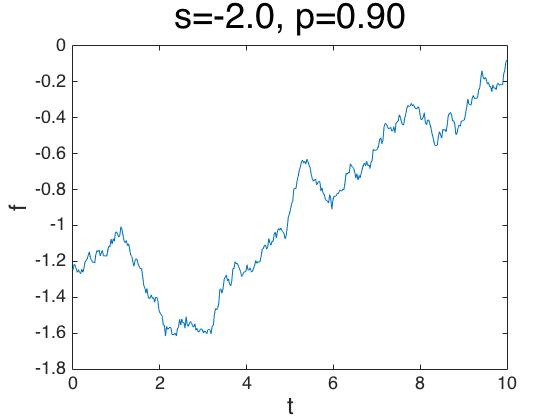} \includegraphics[width=1.12in]{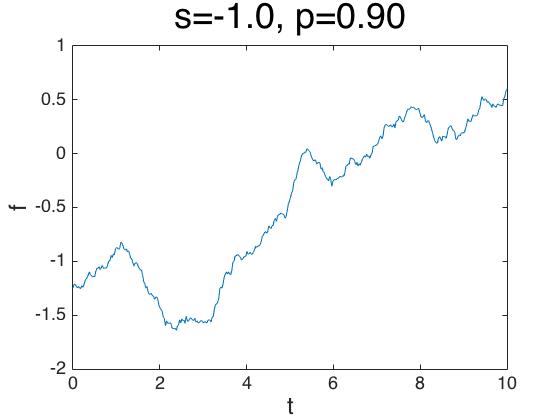} \includegraphics[width=1.12in]{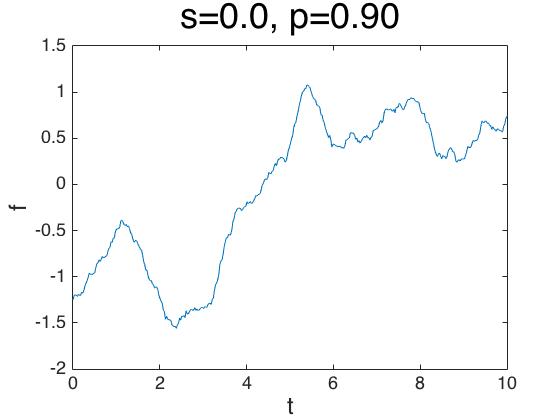} \includegraphics[width=1.12in]{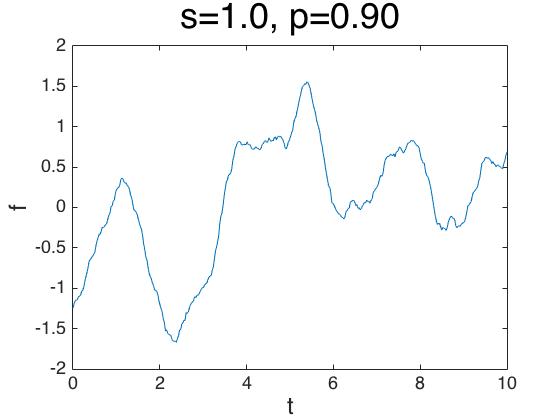} \includegraphics[width=1.12in]{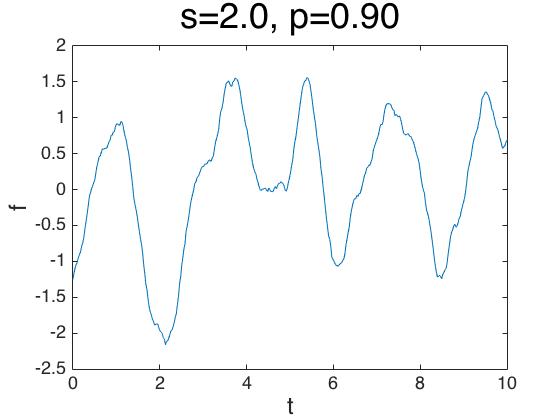} \\
   \includegraphics[width=1.12in]{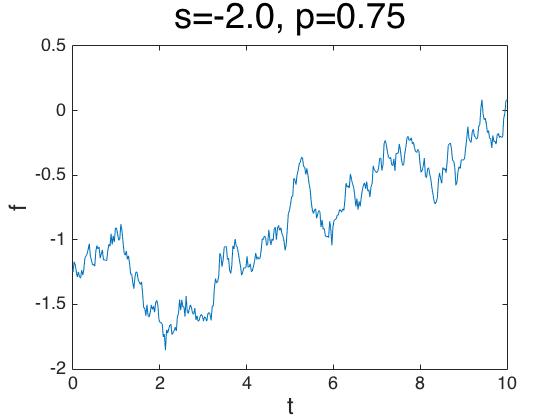} \includegraphics[width=1.12in]{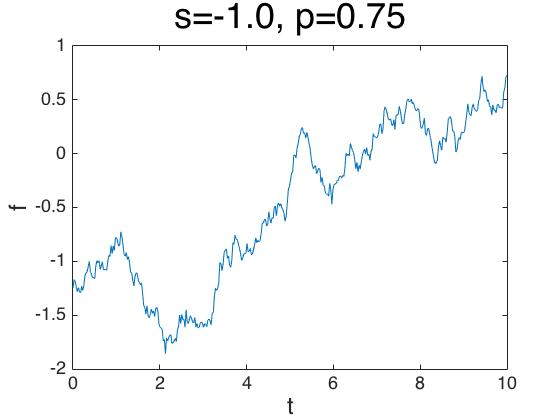} \includegraphics[width=1.12in]{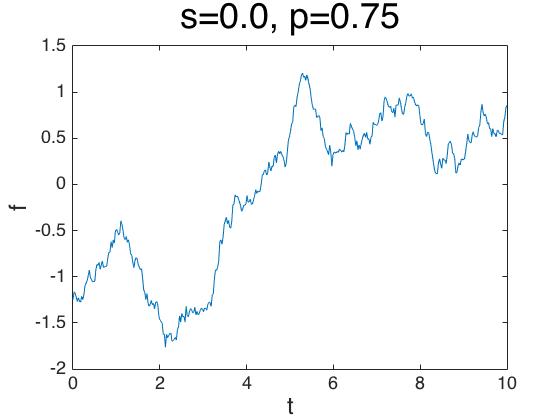} \includegraphics[width=1.12in]{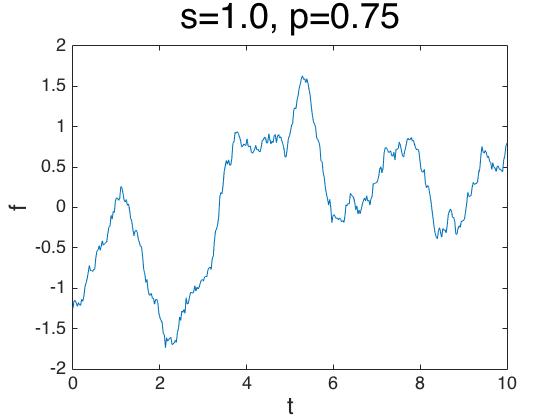} \includegraphics[width=1.12in]{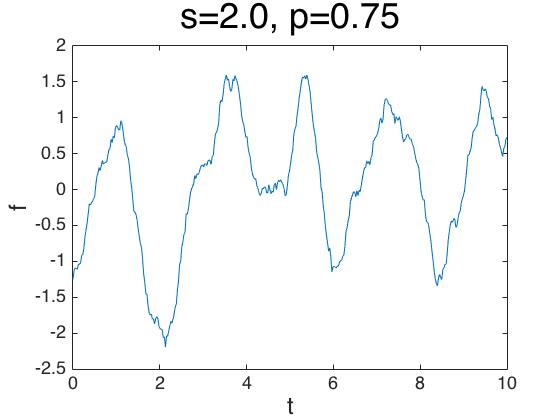} \\
   \includegraphics[width=1.12in]{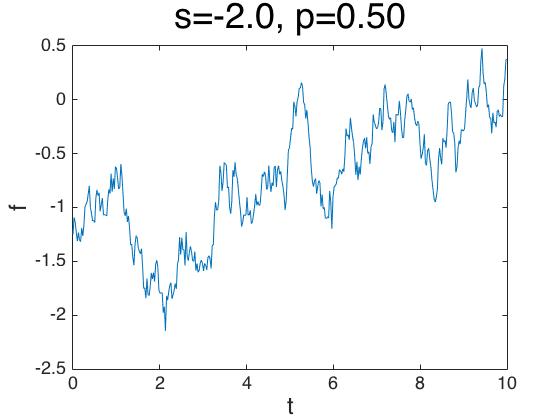} \includegraphics[width=1.12in]{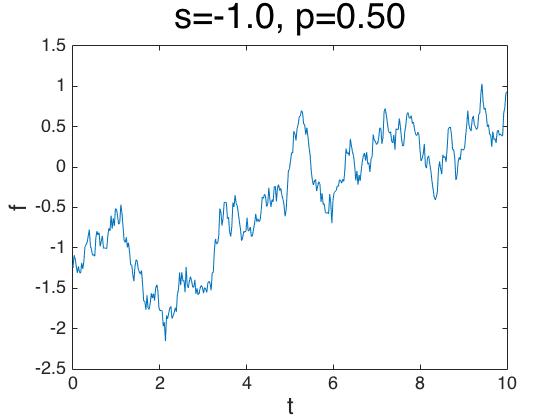} \includegraphics[width=1.12in]{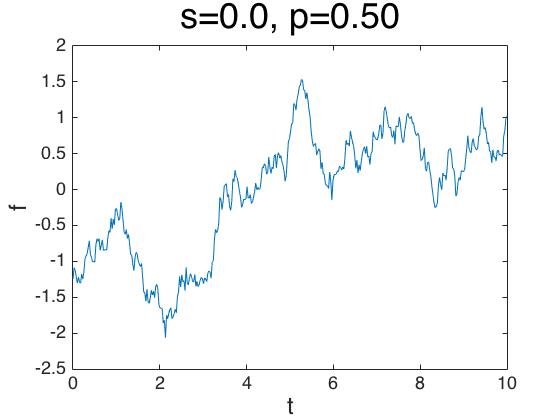} \includegraphics[width=1.12in]{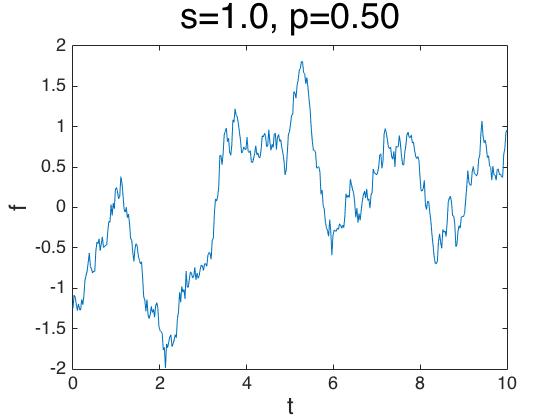} \includegraphics[width=1.12in]{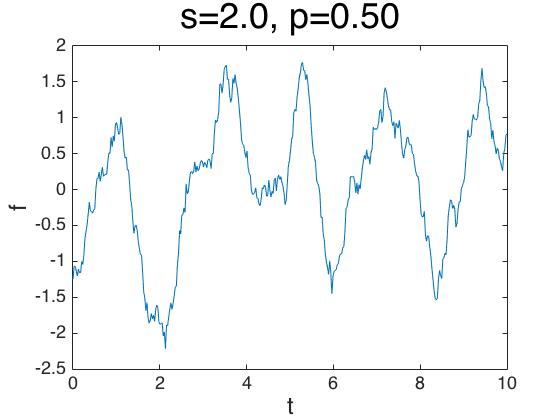} \\
   \includegraphics[width=1.12in]{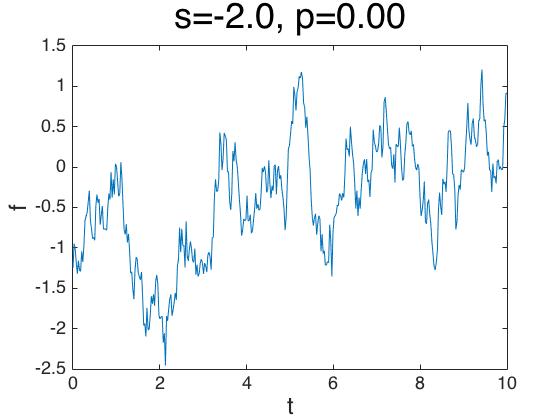} \includegraphics[width=1.12in]{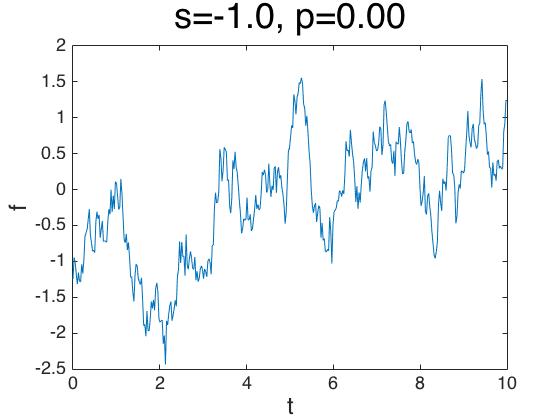} \includegraphics[width=1.12in]{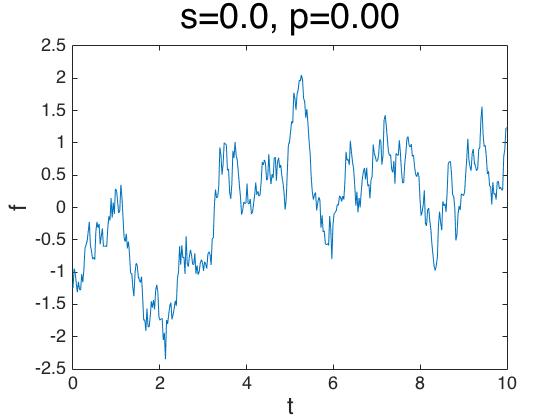} \includegraphics[width=1.12in]{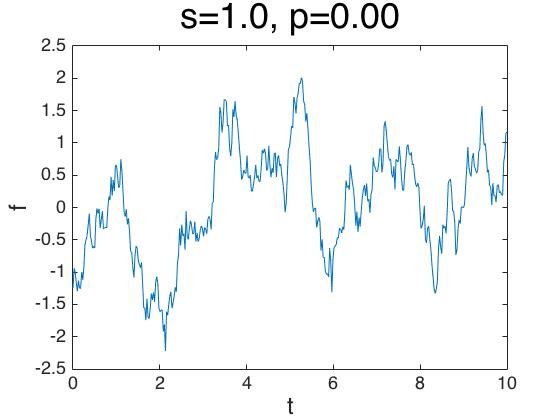} \includegraphics[width=1.12in]{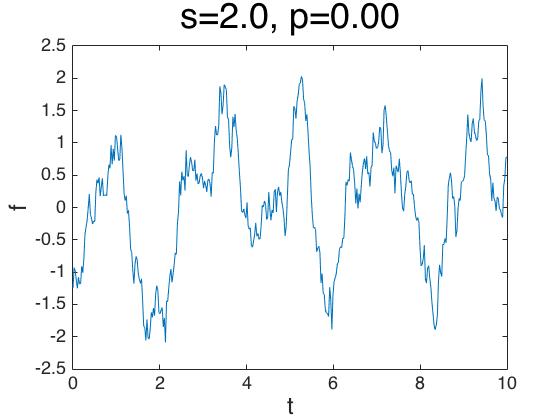} \\
   \includegraphics[width=1.12in]{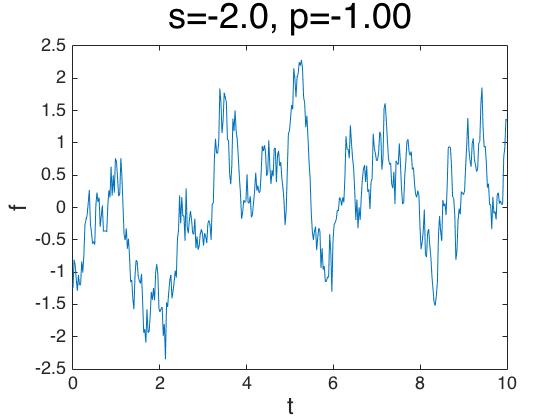} \includegraphics[width=1.12in]{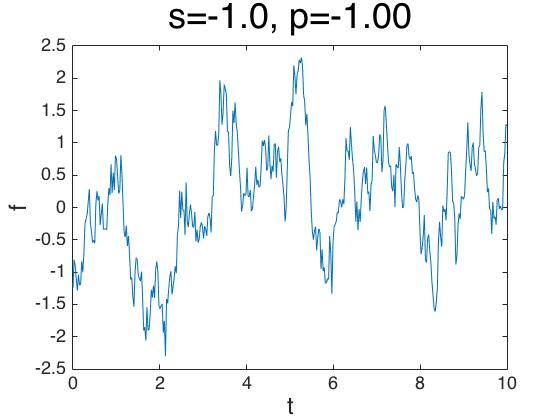} \includegraphics[width=1.12in]{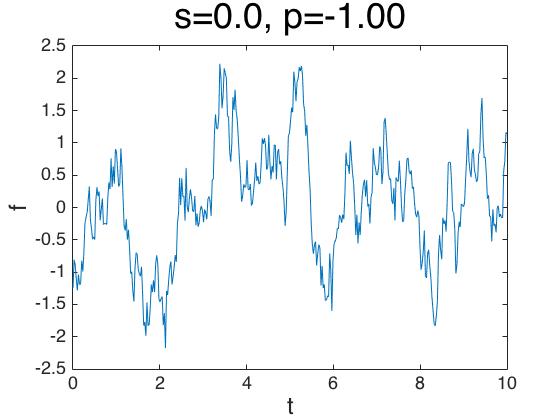} \includegraphics[width=1.12in]{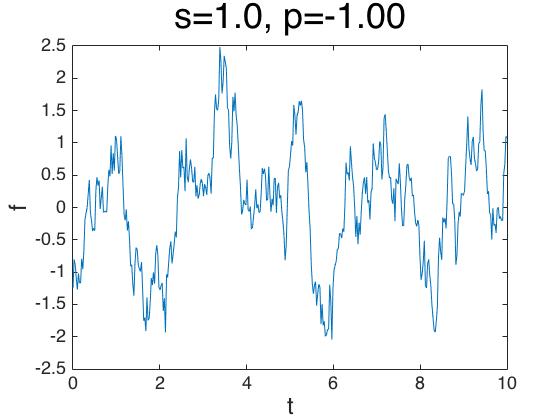} \includegraphics[width=1.12in]{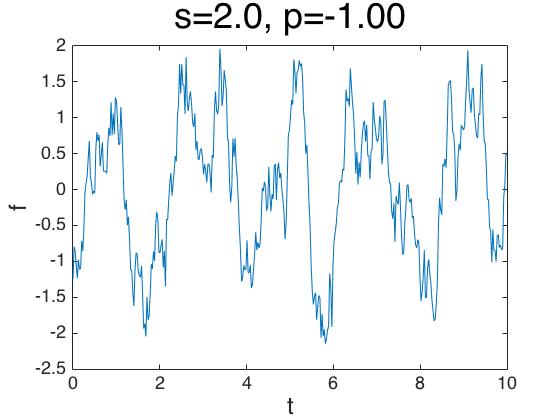} \\
   \caption{Samples from the 2Dsys family on $t=[0,10]$ for $k=h=0$ with $s=-2,-1,0,+1,+2$ horizontally, $p=-1,0,0.5,0.75,0.9,1$ vertically and the same random seed in each case.}
   \label{fig:samples}
\end{figure}

Figure~\ref{fig:samples} shows some samples generated by the family.  We fixed $k=h=0$ as their effects can be removed by changing the scales.  The samples were generated using the GPML script in the Supplementary Information.

Together, the parameters $j$ and $s$ are able to control both the roughness of the signal as well as the coherence of the oscillations, and the family includes regimes of no oscillations, low Q-factor oscillations and high Q-factor oscillations (for a second-order linear system, the {\em quality factor} $Q= \sqrt{\det}/(A+D) = \frac12 e^{s/2}$). Thus, the proposed covariance function can capture a diverse range of behaviour with only 4 parameters.

A special case of our family is OUosc of \cite{PMPR}:~it has $j=0$, $\Delta < 0$ ($s>0$).  Other special cases appear in \cite{M}.

\section{Applications}
The covariance function 2Dsys can be used to fit data to solutions of the system of two stochastic differential equations if parameters are known or assumed.  It can also be used to fit parameters of the model.

For typical applications, one might also want to add a mean, making a five-parameter GP.

One significant application of the covariance function is to deciding whether a time-series corresponds to forced damped oscillation or a forced overdamped system.  This could be done by examining the posterior probability over the parameters.  The part of parameter space with $s>0$ corresponds to oscillatory response, the part with $s<0$ to non-oscillatory response.  So it suffices to compute the posterior probability of $s>0$.  The odds on oscillatory response are given by
\begin{equation}
\frac{P(osc|D)}{P(not|D)} = \frac{\int_{s>0} P(D|h,s,k,j) d\pi(h,s,k,j)}{\int_{s<0} P(D|h,s,k,j) d\pi(h,s,k,j)},
\end{equation}
where data $D$ is a sequence of $n$ values $x_i$ at times $t_i$, 
\begin{equation}
P(D|h,s,k,j) = \frac{\exp(-\frac12 x^T c^{-1} x)}{(2\pi)^{n/2} (\det c)^{1/2}},
\end{equation}
with $c_{ij} = C(t_i,t_j)$, and $d\pi$ is a prior probability distribution on the parameters.

\section{More}
The 2D linear stochastic system defines more generally a GP on $\R \times\{1,2\}$, representing the vector-valued process in $\R^2$.
In contexts where two observables are measured it makes sense to use this more general family of covariance functions, given by (\ref{eq:C}), but it has 7 parameters (and one should add two means).

More generally, a better model might be a stochastic linear system with many components.  Inferring modes of oscillation from such a model is the topic of \cite{M}.

\section*{Acknowledgements}
We thank Magnus Rattray for putting us in touch with each other.
RM is grateful for the support of the Alan Turing Institute under a Fellowship award TU/B/000101.

\section*{Supplementary Information}
\noindent GPML script for the family of 2Dsys covariance functions.

\vskip 1ex
\noindent {\texttt
function [K,dK] = cov2D(hyp, x, z)

%
%
%
%
%

if nargin$<$2, K = '4'; return; end                  

if nargin$<$3, z = []; end                                   

xeqz = isempty(z); dg = strcmp(z,'diag');                       

[n,D] = size(x);

if D$>$1, error('Covariance is defined for 1d input only.'), end

sigma   = exp(hyp(1));

Delta = sigma$^2$*(1-exp(hyp(2))); sq=sqrt(Delta);

S11 = exp(2*hyp(3));

p=hyp(4); j = sin(pi*p/2);

if dg,  T = zeros(size(x,1),1);
  
else

\quad  if xeqz,   T = bsxfun(@plus,x,-x');
    
\quad  else   T = bsxfun(@plus,x,-z');
    
\quad  end
  
end

if Delta,   K = S11*exp(-sigma*abs(T)).*(cosh(sq*T) + sigma*j*sinh(sq*abs(T))/sq); 
    
else K = S11*exp(-sigma*abs(T)).*(1.0+j*sigma*abs(T));

end

if nargout $>$ 1,  dK = @(Q) dirder(Q,K,T,x,sigma,sq,S11,j,p,dg,xeqz);
  
\noindent end

\vskip 1ex

\noindent function [dhyp,dx] = dirder(Q,K,T,x,mu,sq,a2,j,p,dg,xeqz)

  theta = pi*p/2;
  
  A=-mu*abs(T).*K+a2*exp(-mu*abs(T)).*(sq*T.*sinh(sq*T)+mu*j*abs(T).*cosh(sq*T));
  
  B=(sq$^2$-mu$^2$)*mu/(2*sq)*a2*exp(-mu*abs(T)).*(T.*sinh(sq*T)+mu*j*abs(T).*cosh(sq*T)/sq-mu*j*sinh(sq*abs(T))/sq$^2$);
  
  C=a2*exp(-mu*abs(T)).*mu*cos(theta)*sinh(sq*abs(T))*pi/(2*sq);
  
  dhyp = [A(:)'*Q(:);   B(:)'*Q(:);  2*K(:)'*Q(:);  C(:)'*Q(:)];
          
  if nargout $>$ 1
  
\quad    R = a2*exp(-mu*abs(T)).*((sq-mu$^2$*j/sq)*sinh(sq*T)+mu*(j-1)*sign(T).*cosh(sq*T)).*Q;
    
    if dg,   dx = zeros(size(x));
      
    else
    
\quad     if xeqz,   dx = sum(R,2)-sum(R,1)';
        
\quad     else   dx = sum(R,2);
        
\quad      end
      
    end
    
\noindent  end
}

\end{document}